\documentclass[a4paper,11pt]{amsart}
\usepackage{graphicx}
\usepackage{mathptmx}
\usepackage{amsmath}
\usepackage{amssymb}
\usepackage{enumitem}
\usepackage{xcolor}
\usepackage{xparse}
\NewDocumentCommand{\eulerian}{omm}
 {%
  \genfrac<>{0pt}{}{#2}{#3}%
  \IfValueT{#1}{_{\!#1}}%
 }

 \newmuskip\pFqmuskip

\newcommand*\pFq[6][8]{%
  \begingroup 
  \pFqmuskip=#1mu\relax
  \mathchardef\normalcomma=\mathcode`,
  \mathcode`\,=\string"8000
  \begingroup\lccode`\~=`\,
  \lowercase{\endgroup\let~}\pFqcomma
  {}_{#2}F_{#3}{\left(\genfrac..{0pt}{}{#4}{#5}\bigg|#6\right)}%
  \endgroup
}
\newcommand{\pFqcomma}{{\normalcomma}\mskip\pFqmuskip}

\newtheorem{theorem}{Theorem}

\newtheorem{corollary}[theorem]{Corollary}

\begin{document}

\title[Some identities related to degenerate $r$-Bell and degenerate Fubini polynomials]{Some identities related to degenerate $r$-Bell and degenerate Fubini polynomials}

\author{Taekyun  Kim}
\address{Department of Mathematics, Kwangwoon University, Seoul 139-701, Republic of Korea}
\email{tkkim@kw.ac.kr}

\author{Dae San Kim}
\address{Department of Mathematics, Sogang University, Seoul 121-742, Republic of Korea}
\email{dskim@sogang.ac.kr}

\author{Jongkyum Kwon}
\address{Department of Mathematics Education, Gyeongsang National University, Jinju 52828, Republic of Korea}
\email{mathkjk26@gnu.ac.kr}



\begin{abstract}
Many works have been done in recent years as to explorations for degenerate versions of some special polynomials and numbers, which began with the introduction of the degenerate Bernoulli and degenerate Euler polynomials by Carlitz. The aim of this paper is to study some properties, recurrence relations and identities related to the degenerate $r$-Bell polynomials, the two variable degenerate Fubini polynomials and the degenerate $r$-Stirling numbers of the second kind. Especially, we express the power series $\sum_{n=0}^{\infty}\sum_{k=0}^{n}(k+r)_{p,\lambda}\frac{x^n}{n!}$ in terms of the degenerate $r$-Bell polynomials, of the degenerate $r$-Stirling numbers of the second kind and of the degenerate Fubini polynomials.
\end{abstract}

\maketitle

\section{Introduction}
It is Carlitz who began to explore degenerate versions of Bernoulli and Euler polynomials, namely the degenerate Bernoulli and degenerate Euler polynomials (see [2]). These explorations for degenerate versions of some polynomials and numbers have regained interests of some mathematicians in recent years and yielded many interesting results by using various tools in their study. These tools include the generating functions, $p$-adic analysis, umbral calculus, combinatorial methods, probability theory, differential equations, analytic number theory and operator theory (see [5, 7-12]).\par
The aim of this paper is to study by means of generating functions some properties, recurrence relations and identities related to the degenerate $r$-Bell polynomials, the two variable degenerate Fubini polynomials and the degenerate $r$-Stirling numbers of the second kind. \par
In more detail, the outline of this paper is as follows. In Section 1, we recall the degenerate exponential functions $e_{\lambda}^{x}(t)$ in \eqref{1}, the degenerate $r$-Stirling numbers of the second kind ${n+r \brace k+r}_{r,\lambda}$ in \eqref{3}, the degenerate $r$-Bell polynomials $\phi_{n,\lambda}^{(r)}(x)$ in \eqref{5} and the two variable degenerate Fubini polynomials $F_{n,\lambda}(x|r)$ in \eqref{10}. The Section 2 is the main results of this paper. We derive recurrence relations for the degenerate $r$-Bell polynomilas in \eqref{16} and Theorem 1 by considering the derivative of the generating function of those polynomials. In Corollary 3, which is equivalent to Theorem 2, we show that the antiderivative of $\phi_{p,\lambda}^{(r)}(x)$ is given by $e^{-x}\sum_{n=1}^{\infty}\Big((r)_{p,\lambda}+(1+r)_{p,\lambda}+\cdots+(n-1+r)_{p,\lambda}\Big)\frac{x^{n}}{n!}$. We deduce three different expressions for $\sum_{n=0}^{\infty}\Big((r)_{p,\lambda}+(1+r)_{p,\lambda}+\cdots+(n+r)_{p,\lambda}\Big)\frac{x^{n}}{n!}$. Indeed, we represent this power series in terms of the degenerate $r$-Bell polynomials in \eqref{31}, of the degenerate $r$-Stirling numbers of the second kind in Theorem 4 and of the degenerate Fubini polynomials in Theorem 6. We also find some identities involving the two variable degenerate Fubini polynomials in Corollary 7, Theorem 8 and Theorem 10.

For any nonzero $\lambda\in\mathbb{R}$, the degenerate exponential functions are defined by Kim-Kim as
\begin{equation}
e_{\lambda}^{x}(t)=\sum_{n=0}^{\infty}\frac{(x)_{k,\lambda}}{k!}t^{n},\quad (\mathrm{see}\ [5,12]),\label{1}
\end{equation}
where the generalized falling factorial sequence is given by
\begin{equation*}
(x)_{0,\lambda}=1,\quad (x)_{n,\lambda}=x(x-\lambda)\cdots (x-(n-1)\lambda),\quad (n\ge 1).
\end{equation*}
From \eqref{1}, we note that
\begin{equation}
e_{\lambda}^{x}(t)=(1+\lambda t)^{\frac{x}{\lambda}},\quad e_{\lambda}(t)=e_{\lambda}^{1}(t)=(1+\lambda t)^{\frac{1}{\lambda}},\quad\mathrm{and}\quad \lim_{\lambda\rightarrow 0}e_{\lambda}^{x}(t)=e^{xt}.\label{2}
\end{equation}
Let $r$ be a nonnegative integer. The degenerate $r$-Stirling numbers of the second kind are defined by means of the degenerate falling factorial sequence as
\begin{equation}
(x+r)_{n,\lambda}=\sum_{k=0}^{n}{n+r \brace k+r}_{r,\lambda}(x)_{k},\quad (n\ge 0),\quad (\mathrm{see}\ [7,9,12]). \label{3}
\end{equation}
It follows from \eqref{3} that
\begin{equation}
\frac{1}{k!}\Big(e_{\lambda}(t)-1\Big)^{k}e_{\lambda}^{r}(t)=\sum_{n=k}^{\infty}{n+r \brace k+r}_{r,\lambda}\frac{t^{n}}{n!},\quad (n\ge 0),\quad (\mathrm{see}\ [7,9,12]). \label{4}
\end{equation}
Note that $\displaystyle \lim_{\lambda\rightarrow 0}{n+r \brace k+r}_{r,\lambda}={n+r \brace k+r}_{r}\displaystyle$ are the $r$-Stirling numbers of the second kind. It is well known that the $r$-Stirling numbers of the second kind ${n+r \brace k+r}_{r}$ counts the number of partitions of the set $[n]=\{1,2,\dots,n\}$ into nonempty subsets in such a way that the numbers $1,2,3,\dots,r$ are in distinct subsets (see [2,3,7,9,12]). \par
In [12], the degenerate $r$-Bell polynomials $\phi_{n,\lambda}^{(r)}(x)$ are given by
\begin{equation}
e_{\lambda}^{r}(t)e^{x(e_{\lambda}(t)-1)}=\sum_{n=0}^{\infty}\phi_{n,\lambda}^{(r)}(x)\frac{t^{n}}{n!}. \label{5}
\end{equation}
From \eqref{4} and \eqref{5}, we note that
\begin{equation}
\phi_{n,\lambda}^{(r)}(x)=\sum_{k=0}^{n}{n+r \brace k+r}_{r,\lambda}x^{k},\quad (n\ge 0),\quad (\mathrm{see}\ [6,7-10,12]).\label{6}
\end{equation}
By \eqref{5}, we get
\begin{equation}
\begin{aligned}
    \sum_{n=0}^{\infty}\phi_{n,\lambda}^{(r)}(x)\frac{t^{n}}{n!} &=e^{-x}\sum_{k=0}^{\infty}x^{k}\frac{1}{k!}e_{\lambda}^{k+r}(t) \\
    &=\sum_{n=0}^{\infty}\bigg(e^{-x}\sum_{k=0}^{\infty}\frac{(k+r)_{n,\lambda}}{k!}x^{k}\bigg)\frac{t^{n}}{n!}.
\end{aligned}\label{7}
\end{equation}
By comparing the coefficients on both sides of \eqref{7}, we get the Dobinski-like formula
\begin{equation}
\phi_{n,\lambda}^{(r)}(x)=e^{-x}\sum_{k=0}^{\infty}\frac{(k+r)_{n,\lambda}}{k!}x^{k},\quad (n\ge 0).\label{8}
\end{equation}
It is known that the two variable Fubini polynomials are defined by
\begin{equation}
\frac{1}{1-x(e^{t}-1)}e^{yt}=\sum_{n=0}^{\infty}F_{n}(x|y)\frac{t^{n}}{n!},\quad (\mathrm{see}\ [1,4,13,15]). \label{9}
\end{equation}
When $y=0$, $F_{n}(x|0)$ are called the ordinary Fubini polynomials (see [5,10,11,14]). \par
Recently, the two variable degenerate Fubini polynomials are introduced as
\begin{equation}
\frac{1}{1-x(e_{\lambda}(t)-1)}e_{\lambda}^{r}(t)=\sum_{n=0}^{\infty}F_{n,\lambda}(x|r)\frac{t^{n}}{n!},\quad (\mathrm{see}\ [5]). \label{10}
\end{equation}
When $r=0$, $F_{n,\lambda}(x|0)=F_{n,\lambda}(x)$ are called the degenerate Fubini polynomials (see [5,8,10,11]). \par
In this paper, we give some new identities of degenerate $r$-Bell polynomials which are derived from the properties of two variable degenerate Fubini polynomials.

\section{Some identities of degenerate $r$-Bell polynomials}
By \eqref{4} and \eqref{10}, we get
\begin{align}
\sum_{n=0}^{\infty}F_{n,\lambda}(x|r)\frac{t^{n}}{n!}&=\frac{1}{1-x(e_{\lambda}(t)-1)}e_{\lambda}^{r}(t)=\sum_{k=0}^{\infty}x^{k}\Big(e_{\lambda}(t)-1\Big)^{k}e_{\lambda}^{r}(t) \label{11}\\
&=\sum_{k=0}^{\infty}x^{k}k!\sum_{n=k}^{\infty}{n+r \brace k+r}_{r,\lambda}\frac{t^{n}}{n!} \nonumber \\
&=\sum_{n=0}^{\infty}\bigg(\sum_{k=0}^{n}x^{k}k!{n+r \brace k+r}_{r,\lambda}\bigg)\frac{t^{n}}{n!}.\nonumber
\end{align}
Comparing the coefficients on both sides of \eqref{11}, we have
\begin{equation}
F_{n,\lambda}(x|r)=\sum_{k=0}^{n}k!{n+r \brace k+r}_{r,\lambda}x^{k},\quad (n\ge 0).\label{12}
\end{equation}
From \eqref{5}, we have
\begin{align}
\frac{d}{dx}e^{x(e_{\lambda}(t)-1)}\cdot e_{\lambda}^{r}(t)&=\Big(e_{\lambda}(t)-1\Big)e_{\lambda}^{r}(t)e^{x(e_{\lambda}(t)-1)} \label{13} \\
&=\Big(e_{\lambda}(t)-1\Big)\sum_{l=0}^{\infty}\phi_{l,\lambda}^{(r)}(x)\frac{t^{l}}{l!}\nonumber \\
&=\sum_{n=0}^{\infty}\bigg(\sum_{k=0}^{n}\binom{n}{k}(1)_{n-k,\lambda}\phi_{k,\lambda}^{(r)}(x)-\phi_{n,\lambda}^{(r)}(x)\bigg)\frac{t^{n}}{n!}. \nonumber
\end{align}
Thus, by \eqref{13}, we get
\begin{equation}
\begin{aligned}
    \frac{d}{dx}\phi_{n,\lambda}^{(r)}(x)&=\sum_{k=0}^{n}\binom{n}{k}(1)_{n-k,\lambda}\phi_{k,\lambda}^{(r)}(x)-\phi_{n,\lambda}(x) \\
    &=\sum_{k=0}^{n-1}\binom{n}{k}(1)_{n-k,\lambda}\phi_{k,\lambda}^{(r)}(x),\quad (n\ge 1).
\end{aligned}   \label{14}
\end{equation}
From \eqref{2} and \eqref{5}, we note that
\begin{align}
&\sum_{n=0}^{\infty}\phi_{n+1,\lambda}^{(r)}(x)\frac{t^{n}}{n!}=\frac{d}{dt}\Big(e^{x(e_{\lambda}(t)-1)}\cdot e_{\lambda}^{r}(t)\Big) \label{15} \\
&\quad =xe_{\lambda}^{1-\lambda}(t)e^{x(e_{\lambda}(t)-1)}e_{\lambda}^{r}(t)+e^{x(e_{\lambda}(t)-1)} re_{\lambda}^{r-\lambda}(t) \nonumber \\
&\quad =xe_{\lambda}^{1-\lambda}(t)e^{x(e_{\lambda}(t)-1)} e_{\lambda}^{r}(t)+e^{x(e_{\lambda}(t)-1)} e_{\lambda}^{r-1}(t)re_{\lambda}^{1-\lambda}(t) \nonumber \\
&\quad =x\sum_{l=0}^{\infty}(1-\lambda)_{l,\lambda}\frac{t^{l}}{l!}\sum_{k=0}^{\infty}\phi_{k,\lambda}^{(r)}(x)\frac{t^{k}}{k!}+r\sum_{l=0}^{\infty}(1-\lambda)_{l,\lambda}\frac{t^{l}}{l!}\sum_{k=0}^{\infty}\phi_{k,\lambda}^{(r-1)}(x)\frac{t^{k}}{k!} \nonumber \\
&\quad=\sum_{n=0}^{\infty}\bigg\{\sum_{k=0}^{n}\binom{n}{k}(1-\lambda)_{n-k,\lambda}(x\phi_{k,\lambda}^{(r)}(x)+r\phi_{k,\lambda}^{(r-1)}(x)\Big)\bigg\}\frac{t^{n}}{n!}.\nonumber
\end{align}
Comparing the coefficients on both sides of \eqref{15}, we have
\begin{equation}
\phi_{n+1,\lambda}^{(r)}(x)=\sum_{k=0}^{n}\binom{n}{k}(1-\lambda)_{n-k,\lambda}\Big(x\phi_{k,\lambda}^{(r)}(x)+r\phi_{k,\lambda}^{(r-1)}(x)\Big),\quad (n\ge 0).\label{16}
\end{equation}
From \eqref{2}, we note that
\begin{equation}
(1-\lambda)_{n-k,\lambda}=(1)_{n-k,\lambda}\big(1-(n-k)\lambda\big). \label{17}
\end{equation}
Thus, by \eqref{16} and \eqref{17}, we get
\begin{align}
\phi_{n+1,\lambda}^{(r)}(x)&=\sum_{k=0}^{n}\binom{n}{k}(1)_{n-k,\lambda}\Big(x\phi_{k,\lambda}^{(r)}(x)+r\phi_{k,\lambda}^{(r-1)}(x)\Big)\big(1-(n-k)\lambda\big)\label{18} \\
&=\sum_{k=0}^{n}\binom{n}{k}(1)_{n-k,\lambda}\Big(x\phi_{k,\lambda}^{(r)}(x)+r\phi_{k,\lambda}^{(r-1)}(x)\Big)-\lambda\sum_{k=0}^{n}\binom{n}{k}(n-k)(1)_{n-k,\lambda}\Big(x\phi_{k,\lambda}^{(r)}(x)+r\phi_{k,\lambda}^{(r-1)}(x)\Big) \nonumber \\
&=\sum_{k=0}^{n}\binom{n}{k}(1)_{n-k,\lambda}\Big(x\phi_{k,\lambda}^{(r)}(x)+r\phi_{k,\lambda}^{(r-1)}(x)\Big)-n\lambda\sum_{k=0}^{n-1}\binom{n-1}{k}(1)_{n-k,\lambda}\Big(x\phi_{k,\lambda}^{(r)}(x)+r\phi_{k,\lambda}^{(r-1)}(x)\Big). \nonumber
\end{align}
Thus, by \eqref{18}, we get
\begin{equation}
\begin{aligned}
    \sum_{k=0}^{n}\binom{n}{k}(1)_{n-k,\lambda}\phi_{k,\lambda}^{(r)}(x)&=\frac{1}{x}\phi_{n+1,\lambda}^{(r)}(x)-\frac{r}{x}\sum_{k=0}^{n}\binom{n}{k}(1)_{n-k,\lambda}\phi_{k,\lambda}^{(r-1)}(x) \\
    &\quad +n\lambda\sum_{k=0}^{n-1}\binom{n-1}{k}(1)_{n-k,\lambda}\phi_{k,\lambda}^{(r)}(x)+\frac{r}{x}n\lambda\sum_{k=0}^{n-1}\binom{n-1}{k}(1)_{n-k,\lambda}\phi_{k,\lambda}^{(r-1)}(x).
\end{aligned}   \label{19}
\end{equation}
From \eqref{13} and \eqref{19}, we have
\begin{align}
&\frac{1}{x}\phi_{n+1,\lambda}^{(r)}(x)=\sum_{k=0}^{n-1}\binom{n}{k}(1)_{n-k,\lambda}\phi_{k,\lambda}^{(r)}(x)+\phi_{n,\lambda}^{(r)}(x)+\frac{r}{x}\sum_{k=0}^{n-1}\binom{n}{k}(1)_{n-k,\lambda}\phi_{k,\lambda}^{(r-1)}(x) \label{20}\\
&+\frac{r}{x}\phi_{n,\lambda}^{(r-1)}(x)-n\lambda\sum_{k=0}^{n-1}\binom{n-1}{k}(1)_{n-k,\lambda}\phi_{k,\lambda}^{(r)}(x)-\frac{r}{x}n\lambda   \sum_{k=0}^{n-1}\binom{n-1}{k}(1)_{n-k,\lambda}\phi_{k,\lambda}^{(r-1)}(x)\nonumber \\
&=\frac{d}{dx}\phi_{n,\lambda}^{(r)}(x)+\phi_{n,\lambda}^{(r)}(x)+\frac{r}{x}\frac{d}{dx}\phi_{n,\lambda}^{(r-1)}(x)+\frac{r}{x}\phi_{n,\lambda}^{(r-1)}(x)-n\lambda\sum_{k=0}^{n-1}\binom{n-1}{k}(1)_{n-k,\lambda}\phi_{k,\lambda}^{(r)}(x)\nonumber \\
&\quad -\frac{r}{x}n\lambda\sum_{k=0}^{n-1}\binom{n-1}{k}(1)_{n-k,\lambda}\phi_{k,\lambda}^{(r-1)}(x). \nonumber
\end{align}
Therefore, by \eqref{20}, we obtain the following theorem.
\begin{theorem}
For $n\in\mathbb{Z}$ with $n\ge 1$, we have
\begin{align*}
    &\frac{1}{x}\phi_{n+1,\lambda}^{(r)}(x)=\frac{d}{dx}\phi_{n,\lambda}^{(r)}(x)+\frac{r}{x}\frac{d}{dx}\phi_{n,\lambda}^{(r-1)}(x)+\phi_{n,\lambda}^{(r)}(x)+\frac{r}{x}\phi_{n,\lambda}^{(r-1)}(x) \\
    &-n\lambda\sum_{k=0}^{n-1}\binom{n-1}{k}(1)_{n-k,\lambda}\phi_{k,\lambda}^{(r)}(x)-\frac{r}{x}n\lambda\sum_{k=0}^{n-1}\binom{n-1}{k}(1)_{n-k,\lambda}\phi_{k,\lambda}^{(r-1)}(x).
\end{align*}
\end{theorem}
For $p\in\mathbb{N}\cup\{0\}$, let
\begin{equation}
t_{\lambda,p}(x,r)=\sum_{n=0}^{\infty}\sum_{k=0}^{n}(k+r)_{p,\lambda}\frac{x^{n}}{n!}. \label{21}
\end{equation}
Then we have
\begin{equation}
t_{\lambda,p}(x,r)=(r)_{p,\lambda}e^{x}+\sum_{n=1}^{\infty}\sum_{k=1}^{n}(k+r)_{p,\lambda}\frac{x^{n}}{n!}.\label{22}
\end{equation}
From \eqref{21}, we note that
\begin{equation}
\begin{aligned}
t_{\lambda,p}(x,r)  &=\sum_{k=0}^{\infty}\bigg(\sum_{n=k}^{\infty}\frac{x^{n}}{n!}\bigg)(k+r)_{p,\lambda} \\
&=\sum_{k=0}^{\infty}(k+r)_{p,\lambda}\bigg(e^{x}-\sum_{l=0}^{k-1}\frac{x^{l}}{l!}\bigg).
\end{aligned}\label{23}
\end{equation}
Let
\begin{equation}
y_{p,\lambda}(x,r)=t_{\lambda,p}(x,r)-e^{x}\phi_{p,\lambda}^{(r)}(x). \label{24}
\end{equation}
Then, by \eqref{8}, \eqref{21} and \eqref{24}, we get
\begin{align}
y_{p,\lambda}(x,r)&=\sum_{n=0}^{\infty}\bigg(\sum_{k=0}^{n}(k+r)_{p,\lambda}\bigg)\frac{x^{n}}{n!}-\sum_{n=0}^{\infty}(n+r)_{p,\lambda}\frac{x^{n}}{n!}\label{25}\\
&=\sum_{n=1}^{\infty}\Big((r)_{p,\lambda}+(1+r)_{p,\lambda}+\cdots+(n-1+r)_{p,\lambda}\Big)\frac{x^{n}}{n!}. \nonumber
\end{align}
For $p\in\mathbb{N}\cup\{0\}$ and from \eqref{25}, we have
\begin{align}
\frac{d}{dx}y_{p,\lambda}(x,r)&=\sum_{n=1}^{\infty}\Big((r)_{p,\lambda}+(1+r)_{p,\lambda}+\cdots+(n-1+r)_{p,\lambda}\Big)\frac{x^{n-1}}{(n-1)!}\nonumber\\
&=\sum_{n=0}^{\infty}\Big((r)_{p,\lambda}+(1+r)_{p,\lambda}+\cdots+(n+r)_{p,\lambda}\Big)\frac{x^{n}}{n!} \label{26} \\
&=\sum_{n=1}^{\infty}\Big((r)_{p,\lambda}+(1+r)_{p,\lambda}+\cdots+(n-1+r)_{p,\lambda}\Big)\frac{x^{n}}{n!}+\sum_{n=0}^{\infty}(n+r)_{p,\lambda}\frac{x^{n}}{n!} \nonumber \\
&=y_{p,\lambda}(x,r)+e^{x}\phi_{p,\lambda}^{(r)}(x). \nonumber
\end{align}
Therefore, by \eqref{20}, we obtain the following theorem.
\begin{theorem}
For $p\in\mathbb{N}\cup\{0\}$, we have
\begin{displaymath}
\frac{d}{dx}y_{p,\lambda}(x,r)=y_{p}(x,r)+e^{x}\phi_{p,\lambda}^{(r)}(x).
\end{displaymath}
\end{theorem}
By Theorem 2, we easily get
\begin{align}
\frac{d}{dx}\Big(e^{-x}y_{p,\lambda}(x,r)\Big)&=-e^{-x}y_{p,\lambda}(x,r)+e^{-x}\frac{d}{dx}y_{p,\lambda}(x,r) \label{27} \\
&=e^{-x}\bigg(\frac{d}{dx}y_{p,\lambda}(x,r)-y_{p,\lambda}(x,r)\bigg) \nonumber \\
&=e^{-x}e^{x} \phi_{p,\lambda}^{(r)}(x)= \phi_{p,\lambda}^{(r)}(x).\nonumber
\end{align}
\begin{corollary}
    For $p\in\mathbb{N}\cup\{0\}$, we have
    \begin{equation}
    \frac{d}{dx}\Big(e^{-x}y_{p,\lambda}(x,r)\Big)=\phi_{p,\lambda}^{(r)}(x).\label{28}
    \end{equation}
\end{corollary}
From \eqref{28}, we note that
\begin{equation}
e^{-x}y_{p,\lambda}(x,r)=\int_{0}^{x}\phi_{p,\lambda}^{(r)}(t)dt \label{29}.
\end{equation}
Thus, by \eqref{29}, we get
\begin{equation}
y_{p,\lambda}(x,r)=e^{x}\int_{0}^{x}\phi_{p,\lambda}^{(r)}(t)dt.\label{30}
\end{equation}
From \eqref{21}, \eqref{24} and \eqref{30}, we note that
\begin{align}
\sum_{n=0}^{\infty}\Big((r)_{p,\lambda}+(1+r)_{p,\lambda}+\cdots+(n+r)_{p,\lambda}\Big)\frac{x^{n}}{n!}
&=y_{p,\lambda}(x,r)+e^{x}\phi_{p,\lambda}^{(r)}(x)\label{31} \\
&=e^{x}\phi_{p,\lambda}^{(r)}(x)+e^{x}\int_{0}^{x}\phi_{p,\lambda}^{(r)}(t)dt. \nonumber
\end{align}
Now, we consider the degenerate differential operator which is give by
\begin{equation}
\bigg(x\frac{d}{dx}+r\bigg)_{n,\lambda}=\bigg(x\frac{d}{dx}+r\bigg)\bigg(x\frac{d}{dx}+r-\lambda)\cdots\bigg(x\frac{d}{dx}+r-(n-1)\lambda\bigg),\label{32}
\end{equation}
where $n$ is a positive integer. \par
For $p\in\mathbb{N}\cup\{0\}$, we have
\begin{equation}
\bigg(x\frac{d}{dx}+r\bigg)_{p,\lambda}e^{x}=\sum_{n=0}^{\infty}\frac{(n+r)_{p,\lambda}}{n!}x^{n}   =e^{x}\phi_{p,\lambda}^{(r)}(x).\label{33}
\end{equation}
From \eqref{24} and \eqref{30}, we note that
\begin{align}
t_{\lambda,p}(x,r)&=e^{x}\phi_{p,\lambda}^{(r)}(x)+e^{x}\int_{0}^{x}\phi_{p,\lambda}^{(r)}(t)dt \label{34} \\
&=e^{x}\phi_{p,\lambda}^{(r)}(x)+e^{x}\sum_{k=0}^{p}{p+r \brace k+r}_{r,\lambda}\frac{x^{k+1}}{k+1} \nonumber \\
&=e^{x}\sum_{k=0}^{p}{p+r \brace k+r}_{r,\lambda}\bigg(x^{k}+\frac{x^{k+1}}{k+1}\bigg).\nonumber
\end{align}
By \eqref{21} and \eqref{34}, we get
\begin{equation}
\sum_{n=0}^{\infty}\Big((r)_{p,\lambda}+(1+r)_{p,\lambda}+\cdots+(n+r)_{p,\lambda}\Big)\frac{x^{n}}{n!}=e^{x}\sum_{k=0}^{p}{p+r \brace k+r}_{r,\lambda}\bigg(x^{k}+\frac{x^{k+1}}{k+1}\bigg). \label{35}
\end{equation}
Therefore, by \eqref{35}, we obtain the following theorem.
\begin{theorem}
For $p\in\mathbb{N}\cup\{0\}$, we have
\begin{displaymath}
\sum_{n=0}^{\infty}\Big((r)_{p,\lambda}+(1+r)_{p,\lambda}+\cdots+(n+r)_{p,\lambda}\Big)\frac{x^{n}}{n!}=e^{x}\sum_{k=0}^{p}{p+r \brace k+r}_{r,\lambda}\bigg(x^{k}+\frac{x^{k+1}}{k+1}\bigg).
\end{displaymath}
\end{theorem}
From \eqref{3}, we note that
\begin{equation}
{n+r+1 \brace k+r}_{r,\lambda}={n+r \brace k+r-1}_{r,\lambda}+(k+r-n\lambda){n+r \brace k+r}_{r,\lambda},\label{36}
\end{equation}
where $n,k\in\mathbb{N}$ with $n\ge k$. \par
From \eqref{22}, \eqref{34} and \eqref{36}, we note that
\begin{align}
&(r)_{p,\lambda}e^{x}+\sum_{n=1}^{\infty}\sum_{k=1}^{n}(k+r)_{p,\lambda}\frac{x^{n}}{n!}=t_{\lambda,p}(x,r)=\sum_{n=0}^{\infty}\Big((r)_{p,\lambda}+\cdots+(n+r)_{p,\lambda}\Big)\frac{x^{n}}{n!} \label{37}\\
&=e^{x}\sum_{k=0}^{p}\bigg(x^{k}+\frac{x^{k+1}}{k+1}\bigg){p+r \brace k+r}_{r,\lambda}=e^{x}\sum_{k=0}^{p}x^{k}{p+r \brace k+r}_{r,\lambda}+e^{x}\sum_{k=0}^{p}{p+r \brace k+r}_{r,\lambda}\frac{x^{k+1}}{k+1} \nonumber \\
&=e^{x}(r)_{p,\lambda}+e^{x}\sum_{k=1}^{p+1}x^{k}{p+r \brace k+r}_{r,\lambda}+e^{x}\sum_{k=1}^{p+1}{p+r \brace k-1+r}_{r,\lambda}\frac{x^{k}}{k} \nonumber \\
&=e^{x}(r)_{p,\lambda}+e^{x}\sum_{k=1}^{p+1}\frac{x^{k}}{k}\bigg((k+r-p\lambda){p+r \brace k+r}_{r,\lambda}+{p+r \brace k+r-1}_{r,\lambda}\bigg) \nonumber \\
&\quad +(p\lambda -r)e^{x}\sum_{k=1}^{p+1}\frac{x^{k}}{k}{p+r \brace k+r}_{r,\lambda} \nonumber \\
&=e^{x}(r)_{p,\lambda}+e^{x}\sum_{k=1}^{p+1}\bigg({p+1+r \brace k+r}_{r,\lambda}+(p\lambda-r){p+r \brace k+r}_{r,\lambda}\bigg)\frac{x^{k}}{k}. \nonumber
\end{align}
By comparing the coefficients on both sides of \eqref{37} and from \eqref{6}, we get
\begin{align}
&\sum_{n=1}^{\infty}\bigg(\sum_{k=1}^{n}(k+r)_{p,\lambda}\bigg)\frac{x^{n}}{n!}=e^{x}\sum_{k=1}^{p+1}\bigg({p+1+r\brace k+r}_{r,\lambda}+(p\lambda-r){p+r \brace k+r}_{r,\lambda}\bigg)\frac{x^{k}}{k} \label{38}\\
&=e^{x}\int_{0}^{x}\bigg(\big(\phi_{p+1,\lambda}^{(r)}(t)-(r)_{p+1,\lambda}\big)+(p\lambda-r)\big(\phi_{p,\lambda}^{(r)}(t)-(r)_{p,\lambda}\big)\bigg)\frac{dt}{t}.\nonumber
\end{align}
Therefore, by \eqref{38}, we obtain the following theorem.
\begin{theorem}
For $p\in\mathbb{N}\cup\{0\}$, we have
\begin{align*}
&\sum_{n=1}^{\infty}\Big((1+r)_{p,\lambda}+(2+r)_{p,\lambda}+\cdots+(n+r)_{p,\lambda}\Big)\frac{x^{n}}{n!}\\
&=e^{x}\int_{0}^{x}\bigg(\big(\phi_{p+1,\lambda}^{(r)}(t)-(r)_{p+1,\lambda}\big)+(p\lambda-r)\big(\phi_{p,\lambda}^{(r)}(t)-(r)_{p,\lambda}\big)\bigg)\frac{dt}{t}.
\end{align*}
\end{theorem}
Now, we observe that
\begin{align}
&\sum_{l=0}^{\infty}(l+r)_{p,\lambda}x^{l}=\sum_{l=0}^{\infty}x^{l}\sum_{k=0}^{p}{p+r \brace k+r}_{r,\lambda}(l)_{k} \label{39} \\
&=\sum_{l=0}^{\infty}x^{l}\sum_{k=0}^{p}{p+r\brace k+r}_{r,\lambda}\binom{l}{k}k!=\sum_{k=0}^{p}{p+r \brace k+r}_{r,\lambda}k!\sum_{l=k}^{\infty}\binom{l}{k}x^{l} \nonumber \\
&=\sum_{k=0}^{p}{p+r \brace k+r}_{r,\lambda}k!x^{k}\sum_{l=0}^{\infty}\binom{l+k}{k}x^{l}=\sum_{k=0}^{p}{p+r\brace k+r}_{r,\lambda}k!x^{k}\bigg(\frac{1}{1-x}\bigg)^{k+1} \nonumber \\
&=\frac{1}{1-x}\sum_{k=0}^{p}{p+r \brace k+r}_{r,\lambda}k!\bigg(\frac{x}{1-x}\bigg)^{k}=\frac{1}{1-x}
F_{p,\lambda}\bigg(\frac{x}{1-x}\bigg|r\bigg). \nonumber
\end{align}
For $p\in\mathbb{N}\cup\{0\}$ and by \eqref{39}, we get
\begin{align}
&\sum_{n=0}^{\infty}\Big((r)_{p,\lambda}+(1+r)_{p,\lambda}+\cdots+(n+r)_{p,\lambda}\Big)x^{n}=\sum_{n=0}^{\infty}\bigg(\sum_{k=0}^{n}(k+r)_{p,\lambda}\bigg)x^{n}\label{40}\\
&=\sum_{k=0}^{\infty}(k+r)_{p,\lambda}\sum_{n=k}^{\infty}x^{n}=\frac{1}{1-x}\sum_{k=0}^{\infty}(k+r)_{p,\lambda}x^{k}\nonumber \\
&=\frac{1}{(1-x)^{2}}F_{p,\lambda}\bigg(\frac{x}{1-x}\bigg|r\bigg). \nonumber
\end{align}
It is easy to show that
\begin{equation}
\bigg(x\frac{d}{dx}+r\bigg)_{p,\lambda}\frac{1}{1-x}=\sum_{n=0}^{\infty}(n+r)_{p,\lambda}x^{n}.\label{41}
\end{equation}
Therefore, by \eqref{40} and \eqref{41}, we obtain the following theorem.
\begin{theorem}
For $p\in\mathbb{N}\cup\{0\}$, we have
\begin{displaymath}
\bigg(x\frac{d}{dx}+r\bigg)_{p,\lambda}\frac{1}{1-x}=\sum_{n=0}^{\infty}(n+r)_{p,\lambda}x^{n}=\frac{1}{1-x}F_{p,\lambda}\bigg(\frac{x}{1-x}\bigg|r\bigg).
\end{displaymath}
In addition, we have
\begin{displaymath}
\sum_{n=0}^{\infty}\Big((r)_{p,\lambda}+(1+r)_{p,\lambda}+\cdots+(n+r)_{p,\lambda}\Big)x^{n}=\frac{1}{(1-x)^{2}}F_{p,\lambda}\bigg(\frac{x}{1-x}\bigg|r\bigg).
\end{displaymath}
\end{theorem}
Now, we observe from Theorem 6 that
\begin{align}
&\sum_{n=0}^{\infty}(n+r)_{p,\lambda}\bigg(\frac{1}{1-x}-\sum_{l=0}^{n}x^{l}\bigg)=\sum_{n=0}^{\infty}(n+r)_{p,\lambda}\frac{x^{n+1}}{1-x} \label{42}\\
&=\frac{x}{1-x}\sum_{n=0}^{\infty}(n+r)_{p,\lambda}x^{n}=\frac{x}{(1-x)^{2}}F_{p,\lambda}\bigg(\frac{x}{1-x}\bigg|r\bigg). \nonumber
\end{align}
\begin{corollary}
For $p\in\mathbb{N}\cup\{0\}$, we have
\begin{displaymath}
\sum_{n=0}^{\infty}(n+r)_{p,\lambda}\bigg(\frac{1}{1-x}-\sum_{l=0}^{n}x^{l}\bigg)= \frac{x}{(1-x)^{2}}F_{p,\lambda}\bigg(\frac{x}{1-x}\bigg|r\bigg).
\end{displaymath}
\end{corollary}
From \eqref{12}, we note that
\begin{align}
&\frac{1}{m!}\bigg(\frac{d}{dx}\bigg)^{m}\bigg(x^{m}F_{n,\lambda}(x|r)\bigg)=\frac{1}{m!}\bigg(\frac{d}{dx}\bigg)^{m}\bigg(\sum_{k=0}^{n}{n+r \brace k+r}_{r,\lambda}k!x^{k+m}\bigg) \label{43} \\
&=\frac{1}{m!}\sum_{k=0}^{n}{n+r \brace k+r}_{r,\lambda}k!(k+m)_{m}x^{k}=\sum_{k=0}^{n}{n+r \brace k+r}_{r,\lambda}k!\binom{k+m}{m}x^{k}.\nonumber
\end{align}
Therefore, by \eqref{43}, we obtain the following theorem.
\begin{theorem}
For $n,m\ge 0$, we have
\begin{displaymath}
\frac{1}{m!}\bigg(\frac{d}{dx}\bigg)^{m}\bigg(x^{m}F_{n,\lambda}(x|r)\bigg)=\sum_{k=0}^{n}{n+r \brace k+r}_{r,\lambda}k!\binom{k+m}{k}x^{k}.
\end{displaymath}
\end{theorem}
For $p\in\mathbb{N}\cup\{0\}$, we have
\begin{align}
&\sum_{n=0}^{\infty}\binom{n+m}{m}(n+r)_{p,\lambda}x^{n}=\sum_{n=0}^{\infty}\binom{n+m}{m}x^{n}\sum_{k=0}^{p}{p+r \brace k+r}_{r,\lambda}(n)_{k}\label{44}\\
&=\sum_{k=0}^{p}{p+r \brace k+r}_{r,\lambda}k!\sum_{n=k}^{\infty}\binom{n+m}{m}\binom{n}{k}x^{n}=\sum_{k=0}^{p}{p+r \brace k+r}_{r,\lambda}k!\sum_{n=k}^{\infty}\binom{k+m}{m}\binom{n+m}{n-k}x^{n}\nonumber \\
&=\sum_{k=0}^{p}{p+r \brace k+r}_{r,\lambda}k!\binom{k+m}{m}x^{k}\sum_{n=0}^{\infty}\binom{n+m+k}{n}x^{n}\nonumber \\
&=\frac{1}{(1-x)^{m+1}}\sum_{k=0}^{p}{p+r \brace k+r}_{r,\lambda}k!\binom{k+m}{m}\bigg(\frac{x}{1-x}\bigg)^{k}. \nonumber
\end{align}
On the other hand, we also have
\begin{align}
\bigg(x\frac{d}{dx}+r\bigg)_{p,\lambda}\bigg(\frac{1}{1-x}\bigg)^{m+1}&=\sum_{n=0}^{\infty}\binom{n+m}{m}\bigg(x\frac{d}{dx}+r\bigg)_{p,\lambda}x^{n}\label{45} \\
&=\sum_{n=0}^{\infty}\binom{n+m}{n}(n+r)_{p,\lambda}x^{n}.\nonumber
\end{align}
Therefore, by \eqref{44} and \eqref{45}, we obtain the following theorem.
\begin{theorem}
For $p\in\mathbb{N}\cup\{0\}$ and $m\in\mathbb{Z}$ with $m\ge 0$, we have
\begin{align*}
\bigg(x\frac{d}{dx}+r\bigg)_{p,\lambda}\bigg(\frac{1}{1-x}\bigg)^{m+1}&=\sum_{n=0}^{\infty}\binom{n+m}{n}(n+r)_{p,\lambda}x^{n} \nonumber \\
&=\frac{1}{(1-x)^{m+1}}\sum_{k=0}^{p}{p+r \brace k+r}_{r,\lambda}k!\binom{k+m}{m}\bigg(\frac{x}{1-x}\bigg)^{k}.
\end{align*}
\end{theorem}
From Theorem 6, we note that
\begin{align}
&\frac{1}{m!}\bigg(\frac{d}{dx}\bigg)^{m}\bigg[\frac{x^{m}}{(1-x)^{2}}F_{p,\lambda}\bigg(\frac{x}{1-x}\bigg|r\bigg)\bigg]    \label{46}\\
&=\sum_{n=0}^{\infty}\Big((r)_{p,\lambda}+(1+r)_{p,\lambda}+\cdots+(n+r)_{p,\lambda}\Big)\frac{1}{m!}\bigg(\frac{d}{dx}\bigg)^{m}x^{n+m} \nonumber \\
&=\sum_{n=0}^{\infty}\Big((r)_{p,\lambda}+(1+r)_{p,\lambda}+\cdots+(n+r)_{p,\lambda}\Big)\binom{n+m}{n}x^{n}.\nonumber
\end{align}
By Theorem 9, \eqref{45} and \eqref{46}, we get
\begin{align}
&\sum_{k=0}^{\infty}(k+r)_{p,\lambda}\bigg(\bigg(\frac{1}{1-x}\bigg)^{m+1}-\sum_{l=0}^{k}\binom{m+l}{l}x^{l}\bigg)=\sum_{k=0}^{\infty}(k+r)_{p,\lambda}\sum_{n=k+1}^{\infty}\binom{n+m}{n}x^{n} \label{47}\\
&=\sum_{n=1}^{\infty}\binom{n+m}{n}x^{n}\sum_{k=0}^{n-1}(k+r)_{p,\lambda}=\sum_{n=1}^{\infty}\binom{n+m}{n}x^{n}\bigg(\sum_{k=0}^{n}(k+r)_{p,\lambda}-(n+r)_{p,\lambda}\bigg) \nonumber \\
&=\sum_{n=0}^{\infty}\binom{n+m}{n}x^{n}\sum_{k=0}^{n}(k+r)_{p,\lambda}-\sum_{n=0}^{\infty}\binom{n+m}{n}x^{n}(n+r)_{p,\lambda} \nonumber \\
&=\frac{1}{m!}\bigg(\frac{d}{dx}\bigg)^{m}\bigg[\frac{x^{m}}{(1-x)^{2}}F_{p,\lambda}\bigg(\frac{x}{1-x}\bigg|r\bigg)\bigg]-\bigg(x\frac{d}{dx}+r\bigg)_{p,\lambda}\bigg(\frac{1}{1-x}\bigg)^{m+1}.\nonumber
\end{align}
Therefore, by \eqref{47}, we obtain the following theorem.
\begin{theorem}
For $m\in\mathbb{Z}$ with $m\ge 0$ and $p\in\mathbb{N}\cup\{0\}$ we have
\begin{align*}
&\frac{1}{m!}\bigg(\frac{d}{dx}\bigg)^{m}\bigg[\frac{x^{m}}{(1-x)^{2}}F_{p,\lambda}\bigg(\frac{x}{1-x}\bigg|r\bigg)\bigg]-\bigg(x\frac{d}{dx}+r\bigg)_{p,\lambda}\bigg(\frac{1}{1-x}\bigg)^{m+1} \\
&=\sum_{k=0}^{\infty}(k+r)_{p,\lambda}\bigg(\bigg(\frac{1}{1-x}\bigg)^{m+1}-\sum_{l=0}^{k}
\binom{m+l}{l}x^{l}\bigg).
\end{align*}
\end{theorem}

\section{Conclusion}
Carlitz initiated a study of degenerate versions of Bernoulli and Euler polynomials. Some mathematicians have investigated various degenerate versions of many special polynomials and numbers with their regained the interests. \par
In this paper, we studied some properties, recurrence relations and identities in connection with the degenerate $r$-Bell polynomials, the two variable degenerate Fubini polynomials and the degenerate $r$-Stirling numbers of the second kind. For example, we showed that the following identities hold true.
\begin{align*}
\sum_{n=0}^{\infty}\Big((r)_{p,\lambda}+(1+r)_{p,\lambda}+\cdots+(n+r)_{p,\lambda}\Big)\frac{x^{n}}{n!}&=e^{x}\phi_{p,\lambda}^{(r)}(x)+e^{x}\int_{0}^{x}\phi_{p,\lambda}^{(r)}(t)dt \\
&=\frac{1}{(1-x)^{2}}F_{p,\lambda}\bigg(\frac{x}{1-x}\bigg|r\bigg)\\
&=e^{x}\sum_{k=0}^{p}{p+r \brace k+r}_{r,\lambda}\bigg(x^{k}+\frac{x^{k+1}}{k+1}\bigg).
\end{align*}
\indent It is one of our future research projects to continue to study various degenerate versions of many special polynomials and numbers and to find their applications to physics, science and engineering as well as to mathematics.


\begin{thebibliography}{9}
\bibitem{1}
Araci, S. \emph{A new class of Bernoulli polynomials attached to polyexponential functions and related identities.} Adv. Stud. Contemp. Math. (Kyungshang) \textbf{31} (2021), no. 2, 195-204.
\bibitem{2}
Carlitz, L. \emph{Degenerate Stirling, Bernoulli and Eulerian numbers.} Utilitas Math. \textbf{15} (1979), 51-88.
\bibitem{3}
Comtet, L. \emph{Advanced combinatorics. The art of finite and infinite expansions.} Revised and enlarged edition. D. Reidel Publishing Co., Dordrecht, 1974. xi+343 pp. ISBN: 90-277-0441-4.
\bibitem{4}
Kilar, N.; Simsek, Y. \emph{Identities and relations for Fubini type numbers and polynomials via generating functions and $p$-adic integral approach.} Publ. Inst. Math. (Beograd) (N.S.) \textbf{106} (120) (2019), 113-123.
\bibitem{5}
Kim, D. S.; Jang, G.-W.; Kwon, H.-I.; Kim, T. \emph{Two variable higher-order degenerate Fubini polynomials.} Proc. Jangjeon Math. Soc. \textbf{21} (2018), no. 1, 5-22.
\bibitem{6}
Kim, T.; Kim, D. S. \emph{Some results on degenerate Fubini and degenerate Bell polynomials.} Applicable Analysis and Discrete Mathematics. \textbf{2022} 2022 OnLine-First (00):35-35.
https://doi.org/10.2298/AADM200310035K
\bibitem{7}
Kim, T.; Kim, D. S. \emph{Degenerate Whitney numbers of first and second kind of Dowling lattices.} Russ. J. Math. Phys. \textbf{29} (2022), no. 3, 358-377.
\bibitem{8}
Kim, T.; Kim, D. S. \emph{On some degenerate differential and degenerate difference operators.} Russ. J. Math. Phys. \textbf{29} (2022), no. 1, 37-46.
\bibitem{9}
Kim, T.; Kim, D. S. \emph{Degenerate $r$-Whitney numbers and degenerate $r$-Dowling polynomials via boson operators.} Adv. in Appl. Math. \textbf{140} (2022), Paper No. 102394, 21 pp.
\bibitem{10}
Kim, T.; Kim, D. S.; Kim, H. K.; Lee, H. \emph{Some properties on degenerate Fubini polynomials.} Appl. Math. Sci. Eng. \textbf{30} (2022), no. 1, 235-248.
\bibitem{11}
Kim, T.; Kim, D. S.; Lee, H.; Kwon, J. \emph{On degenerate generalized Fubini polynomials.} AIMS Math. \textbf{7} (2022), no. 7, 12227-12240.
\bibitem{12}
Kim, T.; Yao, Y.; Kim, D. S.; Jang, G.-W. \emph{Degenerate $r$-Stirling numbers and $r$-Bell polynomials.} Russ. J. Math. Phys. \textbf{25} (2018), no. 1, 44-58.
\bibitem{13}
Muhyi, A.; Araci, S. \emph{A note on $q$-Fubini-Appell polynomials and related properties.} J. Funct. Spaces 2021, Art. ID 3805809, 9 pp.
\bibitem{14}
Roman, S. \emph{The umbral calculus.} Pure and Applied Mathematics, 111. Academic Press, Inc. [Harcourt Brace Jovanovich, Publishers], New York, 1984. x+193 pp. ISBN: 0-12-594380-6
\bibitem{15}
Simsek, Y. \emph{Construction of generalized Leibnitz type numbers and their properties.} Adv. Stud. Contemp. Math. (Kyungshang) \textbf{31} (2021), no. 3, 311-323.
\end{thebibliography}
\end{document}